\documentclass[12pt,leqno]{amsart}
\usepackage{amssymb}
\usepackage{amsmath,amssymb,color}
\newtheorem{theorem}{Theorem}

\newtheorem{lemma}{Lemma}

\numberwithin{equation}{section}

\newcommand{\ep}{\varepsilon}

\newcommand{\ppp}{\partial}

\newcommand{\R}{\mathbb{R}}
\newcommand{\C}{\mathbb{C}}

\newcommand{\www}{\widetilde}

\newcommand{\ooo}{\overline}
\newcommand{\OOO}{\Omega}

\allowdisplaybreaks

\title
[]
{
Inverse coefficient problem for one-dimensional 
evolution equation with 
vanishing initial condition}

\author{$^1$ Oleg Imanuvilov,  
$^2$ Masahiro Yamamoto}

\thanks{
$^1$ Department of Mathematics, Colorado State University, 101 Weber Building, 
Fort Collins CO 80523-1874, USA 
e-mail: {\tt oleg@math.colostate.edu}
\\
$^2$ Graduate School of Mathematical Sciences, The University
of Tokyo, Komaba, Meguro, Tokyo 153-8914, Japan 
e-mail: {\tt myama@ms.u-tokyo.ac.jp}
}

\date{}
\begin{document}
\maketitle

\begin{abstract}{\it
We consider an inverse problem of
determining a coefficient $p(x)$ of an evolution equation
$\sigma\ppp_tu = a(x)\ppp_x^2u - p(x)u$ for 
$0<x<\ell$ and $0<t<T$, where $\sigma \in \C \setminus \{0\}$,
$\ell>0$ and $T>0$ are arbitrarily given.  
Our main result is the uniqueness: by assuming that 
the zeros of initial value $b(x):= u(0,x)$ on $[0, \ell]$ 
is a finite set and each zero is
of order one at most, if two solutions have the same Cauchy data 
at $x=0$ over $(0,T)$ and the same initial value $b(x)$,
then the coefficient $p(x)$ is uniquely determined on $[0,\ell]$.
}
\end{abstract} 
\baselineskip 18pt

\section{Introduction and main result}
This paper is concerned with the following one dimensional linear evolution
equation
\begin{equation}\label{cop}
\sigma\ppp_tu(t,x) = a(x)\ppp_x^2u(t,x) - p(x)u(t,x), \quad 0<t<T, \, 
0<x<\ell,
\end{equation}
with initial condition 
\begin{equation}\label{cop1}
u(0,x)=b(x), \quad 0<x<\ell.
\end{equation}
We consider the following inverse problem: 

{\it Given 
$u(t, 0)$ and $\ppp_xu(t,0)$ for $0<t<T$, 
determine a coefficient $p(x)$ for $0\le x\le \ell$.}
\\

We remark that we can choose a constant $\sigma \in \C \setminus \{ 0\}$
arbitrarily, and $\sigma=1$ and $\sigma=\sqrt{-1}$ provide a parabolic 
equation and a Schr\"odinger equation respectively.

Our main uniqueness result is 

\begin{theorem}\label{Moon}
{\it Let $\sigma\in \Bbb C\setminus \{0\}$,
$p, q \in C^1[0,\ell], a\in C^2[0,\ell], a(x)>0$ on $[0,\ell]$.  Moreover 
we assume that the set of the zeros of 
$b \in C^3[0,\ell]$ is a finite set and each zero is of at most 
order one.
Let $u$ and $\widetilde u \in C^{1,2}([0,T]\times [0,\ell])$ satisfy
$$
\sigma\ppp_tu(t,x) = a(x)\ppp_x^2u(t,x) - p(x)u(t,x), 
$$
$$
\sigma\ppp_t\www u(t,x) = a(x)\ppp_x^2\www u(t,x) - q(x)\www u(t,x), 
\quad 0<t<T,\, 0<x<\ell, 
$$
and
\begin{equation}\label{(1.7)}
u(0,\cdot) = \www u(0,\cdot) = b \quad \mbox{in $(0,\ell)$}.    
\end{equation}
We assume $u(t,0) = \www u(t,0)$ for $0<t<T$.
Then, $p=q$ on $[0,\ell]$ in each case:
\\
Case A: 
\begin{equation}\label{vodolaz}
\ppp_x u(t,0) = \ppp_x\www{u}(t,0) = 0 \quad \mbox{for
$0<t<T$.}
\end{equation}
\\
Case B:
There exists $\ep\in (0,\ell)$ such that $p=q$ on $[0,\ep]$. 
}

We have the same conclusion in the case where we replace
(1.3) by
\begin{equation}\label{(1.77)}
u(T,\cdot) = \www u(T,\cdot) = b \quad \mbox{in $(0,\ell)$}.
\end{equation}
\end{theorem}

Without any assumption on a function $b$, the uniqueness for the inverse 
problem does not hold.  For example, the case $b\equiv 0$ admits that 
a solution $u$ is identically zero, and so we cannot expect the uniqueness.
Theorem \ref{Moon} implies that the uniqueness holds for $b$ given in 
a dense set in the space $C[0,\ell]$.

The determination problems of spatially varying coefficients 
by Cauchy data on a lateral subboundary and spatial data at fixed time
$t_0$, are classified into important and typical inverse problems.
In the case of $\sigma=1$, we refer to the following works;
If spatial data are taken at $t_0 \in (0,T)$, then 
a pioneering work Bukhgeim and Klibanov \cite{BK} proves the uniqueness.
See also Imanuvilov and Yamamoto \cite{IY98}, Isakov \cite{Is1},
Klibanov \cite{Kli}, Yamamoto \cite{Y09} as for the uniqueness and the 
stability for inverse parabolic problems.
The above works discuss general spatial dimensions, but 
have to assume that the spatial data are given 
at $t_0 \in (0,T)$, which means that the conditions 
(\ref{(1.7)}) and (\ref{(1.77)}) can not be considered.  Under additional conditions, 
we can discuss the cases $t_0=0$ and $t_0=T$ (Imanuvilov and Yamamoto 
\cite{IY230}, \cite{IY231}).
In Imanuvilov and Yamamoto \cite{I-M1} and \cite{I-M2}, the uniqueness is 
proved for the cases $t_0=0$ and $t_0=T$ as long as the spatial dimension is 
one.   On the other hand, all the above works requires that the spatial data
$b$ are strictly positive on $[0, \ell]$ if we want to prove the uniqueness of 
coefficients on $[0,\ell]$.  The uniqueness is still an open problem when 
$b$ vanishes in $(0,\ell)$.  Our main result solves this open problem 
in the one-dimensional case.

The proof is based on the transformation operator (e.g., 
Levitan \cite{Le}) and modifies the argument in \cite{I-M1}.

As for one-dimensional inverse parabolic coefficient problems 
requiring the full boundary condition at $x=0$ and $x=\ell$, we can refer to
Murayama \cite{Mu}, Pierce \cite{Pi}, Suzuki \cite{Su},
Suzuki and Murayama \cite{SM}.
 
\section{ Proof of Theorem \ref{Moon}.}\label{SNOM}
Henceforth we set $\OOO := \{ (x,y);\, 0<y<x<\ell \}$ and 
$\OOO_{x_0}:= \{ (x,y);\, 0<y<x<x_0 \}$ with $x_0 \in (0,\ell)$.

First we show
\\
\begin{lemma}\label{Asol}
{\it 
Let $p, q \in C^1[0,\ell], a\in C^2[0,\ell], a(x)>0$ on $[0,\ell]$.
Then there exists a unique ${ K} \in C^2(\ooo{\OOO})$ 
to the following problem:
\begin{equation}\label{(3.1)}
\left\{ \begin{array}{rl}
& a(x)\ppp_x^2{ K}(x,y) -\ppp_y^2( a(y){ K}(x,y) )= { K}(x,y) (q(x)-p(y)), 
\quad (x,y) \in \OOO, \\
& a(0)\ppp_y {K}(x,0)+a'(0){K}(x,0) = 0, \quad 0<x<\ell, \\
& 2a(x)\frac{d}{dx}{K}(x,x)+a'(x){K}(x,x) = q(x) - p(x),\, K(0,0)=0 \quad 
0<x<\ell.
\end{array}\right.
\end{equation}
}
\end{lemma}

This is a Goursat problem and the proof is standard by means of the 
characteristics.  For example, one can find the proof in Suzuki \cite{Su}.

In view of $ K$, we define an operator $K: L^2(0,\ell) \longrightarrow 
L^2(0,\ell)$, denoted by the same letter, by 
\begin{equation}\label{spartak}
(Kv)(x) := \int^x_0  K(x,y)v(y) dy\quad 0<x<\ell.
\end{equation}

Then, we can prove 
\begin{lemma}\label{main} 
Let $p, q \in C^1[0,\ell], a\in C^2[0,\ell], a(x)>0$ on $[0,\ell]$.
Let $u \in  C^{1,2}([0,T]\times [0,\ell])$ satisfy
\begin{equation}
\sigma\ppp_t{u} - a(x)\ppp_x^2{u} + p(x){u} = 0\quad\mbox{ in}\quad  
(0,T)\times (0,\ell).
\end{equation} 
Then $\www{v}$ given by  
\begin{equation}\label{mk}
\www{v}={u} + Ku
\end{equation}
satisfies 
\begin{equation}\label{Mmk1}
\begin{split}&\sigma\ppp_t\www{v} - a(x)\ppp_x^2\www{v} + q(x)\www{v} = -a(0)
K(x,0)\partial_x u(t,0), \quad 0<t<T, \, 0<x<\ell,\\ 
&\partial_x\www{v}(t,0) = \partial_x{u}(t,0), 
\quad \www{v}(t,0)={u}(t,0), \quad 0<t<T.
\end{split}
\end{equation}
\end{lemma}

This is a transformation operator (e.g., 
Levitan \cite{Le}), and is used for one-dimensional inverse problems 
(e.g., \cite{Su}, \cite{SM}).
However, unlike \cite{Su} and \cite{SM}, our method 
never uses eigenfunction expansions or any spectral properties of full
solutions, and so is applicable even if we do not assume the 
boundary conditions at $x=0$ and $x=\ell$.

{\bf Proof.}
We directly verify that $\www{v}$ given by (\ref{mk}), satisfies 
(\ref{Mmk1}). 
First we differentiate equation (\ref{mk}) with respect to the variable $t:$
\begin{equation}\label{(3.2)}
\ppp_t\www{v}(t,x) = \ppp_t{u}(t,x) + \int^x_0 {K}(x,y)\ppp_t{u}(t,y) dy,
\quad 0<x<\ell, \, 0<t<T.                   
\end{equation}
Next, differentiating (\ref{(3.1)}) with respect to variable $x$,
we have
$$
\ppp_x\www{v}(t,x) = \ppp_x{u}(t,x) +{K}(x,x){u}(t,x) 
+ \int^x_0 \ppp_x{K}(x,y){u}(t,y) dy,
$$
and so
\begin{equation}\label{(3.33)}
\begin{split}
& a(x)\ppp^2_x\int^x_0 {K}(x,y){u}(t,y) dy = 
 a(x)\frac{d}{dx}({K}(x,x)){u}(t,x) + a(x){K}(x,x)\ppp_x{u}(t,x)\\
+& a(x)(\ppp_x{K})(x,x){u}(t,x) + a(x)\int^x_0 \ppp_x^2K(x,y){u}(t,y) dy\\
= & a(x) \frac{d}{dx}({K}(x,x)){u}(t,x) +a(x) {K}(x,x)\ppp_x{u}(t,x)\\
+& a(x) (\ppp_x{K})(x,x){u}(t,x) + \int^x_0 \ppp_y^2(a(y){K}(x,y)){u}(t,y) dy\\
+ & \int^x_0 (q(x){K}(x,y)-{K}(x,y)p(y)){u}(t,y) dy, \quad 0<x<\ell, \, 0<t<T.
\end{split}
\end{equation}
Hence, in view of
$\ppp_y(a(y){K}(x, y))\vert_{y=0} =  0$ for $0<x<\ell$, the integration by 
parts yields
\begin{equation}\label{(3.3)}
\begin{split}
&\int^x_0 \ppp_y^2(a(y){K}(x,y)){u}(t,y) dy \\
= & \left[ (\ppp_y(a(y){K}(x,y)){u}(t,y)\right]^{y=x}_{y=0}
- \int^x_0 \ppp_y(a(y){K}(x,y))\ppp_y{u}(t,y) dy\\
=& (a(x)\ppp_y{K}(x,x)+a'(x){K}(x,x)){u}(t,x) 
- \left[ a(y){K}(x,y)\ppp_y{u}(t,y)\right]^{y=x}_{y=0}\\
+ & \int^x_0 a(y){K}(x,y)\ppp_y^2{u}(t,y) dy \\
=& (a(x)\ppp_y{K}(x,x)+a'(x){K}(x,x)){u}(t,x) -a(x){K}(x,x)\ppp_x{u}(t,x)\\ 
+ & a(0){K}(x,0)\ppp_x{u}(t,0)
+ \int^x_0 a(y) {K}(x,y)\ppp_y^2{u}(t,y) dy.
\end{split}
\end{equation}
By (\ref{(3.33)}) and (\ref{(3.3)}), we have
\begin{equation}\nonumber
\begin{split}
& a(x)\ppp^2_x\int^x_0 {K}(x,y){u}(t,y) dy =
  a(x) \frac{d}{dx}({K}(x,x)){u}(t,x) \\
+& a(x) (\ppp_x{K})(x,x){u}(t,x) 
+  \int^x_0 (q(x){K}(x,y)-{K}(x,y)p(y)){u}(t,y) dy\\
+ &(a(x)\ppp_y{K}(x,x)+a'(x){K}(x,x)){u}(t,x) + a(0){K}(x,0)\ppp_x{u}(t,0)\\
+ & \int^x_0 a(y){K}(x,y)\ppp_y^2{u}(t,y) dy.
\end{split}
\end{equation}
Since
$$
\frac{d}{dx}({K}(x,x)) = (\ppp_x{K})(x,x) + (\ppp_y{K})(x,x),
$$
we can rewrite the above equality as
\begin{equation}\nonumber
\begin{split}
& a(x)\ppp^2_x\int^x_0 {K}(x,y){u}(t,y) dy =
  2a(x) \frac{d}{dx}({K}(x,x)){u}(t,x)+a'(x){K}(x,x)){u}(t,x)  \\
+ & \int^x_0 (q(x){K}(x,y)-{K}(x,y)p(y)){u}(t,y) dy-a(0){K}(x,0)\ppp_x{u}(t,0)
                    \\
+ &\int^x_0 a(y){K}(x,y)\ppp_y^2{u}(t,y) dy.
\end{split}
\end{equation}
Finally, using the boundary condition (\ref{(3.1)}), we obtain
\begin{equation}\label{saboromond}
\begin{split}
& a(x)\ppp^2_x\int^x_0 {K}(x,y){u}(t,y) dy =
  q(x){u}(t,x)-p(x){u}(t,x)  \\
+ & \int^x_0 K(x,y)(q(x)-p(y))u(t,y) dy- a(0){K}(x,0)\ppp_x{u}(t,0)
                    \\
+ & \int^x_0 a(y){K}(x,y)\ppp_y^2{u}(t,y) dy.
\end{split}
\end{equation}

Therefore, since $\sigma\ppp_t{u} - a(x)\ppp_x^2{u} + p(x){u} = 0$ 
in $(0,T)\times (0,\ell)$, we obtain
\begin{align*}
& \sigma \ppp_t\www{v} - a(x)\ppp_x^2\www{v} + q(x)\www{v}\\
=& \sigma \ppp_t{u} + \int^x_0 {K}(x,y)\sigma\ppp_t{u}(t,y) dy
- a(x)\ppp_x^2{u} \\
-& \int^x_0 K(x,y)(q(x)-p(y))u(t,y) dy + \int^x_0 a(y){K}(x,y)
\ppp_y^2{u}(t,y) dy-q(x){u}+p(x){u}               \\
+& q(x){{u}}+q(x)\int^x_0 {K}(x,y){u}(t,y) dy- a(0){K}(x,0)\ppp_x{u}(t,0)
= -a(0){K}(x,0)\ppp_x{u}(t,0).
\end{align*}
The proof of Lemma \ref{main}  is complete.
$\blacksquare$

Now we complete the proof of Theorem \ref{Moon}.
First, in terms of (\ref{vodolaz}) and (\ref{Mmk1}),
the function $\www{v}$ given by 
$$
\www{v}={{u}}+ Ku
$$
satisfies 
\begin{equation}\label{mk1}
\sigma\ppp_t\www{v} - a\ppp_x^2\www{v} + q\www{v} = 0\,\mbox{ in}\, 
(0,T)\times (0,\ell),\quad \partial_x\www{v}(t,0)=0,\, \www{v}(t,0)={u}(t,0), 
\quad 0\le t \le T.
\end{equation}
By (\ref{cop}) and (\ref{mk1}), the function $w:= \www {u}- \www{v}$ satisfies 
\begin{equation}\label{park}
\sigma\partial_t w + a\partial^2_xw + qw=0 \,\mbox{ in}\,\, 
(0,T)\times (0,\ell),\quad w(t,0)=\partial_{x}w(t,0)=0, \quad 
0\le t \le T.
\end{equation}
We apply Lemma \ref{KSOKOL2}  which is shown in Section \ref{SK2} to (\ref{park})
so that 
$$
\www{v}=\www u\quad \mbox{on}\quad (0,T)\times (0,\ell).
$$
This equality and (\ref{mk}) imply 
\begin{equation}\label{main1}
\int_0^xK(x,y){u}(0,y)dy=0\quad \mbox{on}\quad  [0,\ell].
\end{equation} 
We consider only the case (\ref{(1.7)}).  For the case (\ref{(1.77)}),
we can derive 
\begin{equation}\label{Amain1}
\int_0^xK(x,y){u}(T,y)dy=0\quad \mbox{on}\quad  [0,\ell],
\end{equation}
and the proof is similar.

Now we reduce Case A to Case B. 
First assume that $\vert b(x)\vert = \vert u(0,x)\vert \ne 0$ for 
$0 \le x < \ep$: some constant, then we can apply \cite{I-M1} to obtain
$p(x) = q(x)$ for $0<x<\ep$.  Therefore, in this situation  Case A reduced to the  Case B.
Then, we can assume that $b(0) = 0$.  By the assumption, 
$0$ is a zero of $b$ of the first order. 
Hence $b(x)= x \www b(x)$, where $\www b\in C^2[0,\ell]$ and
$\vert\www b(x)\vert \ge \kappa>0$ in some neighborhood of  $0$ with some 
constant $\kappa > 0$.
Using (\ref{saboromond}) and (\ref{park}), we have
\begin{equation}\label{saboromond1}
\begin{split}
&0= a(x)\ppp^2_x\int^x_0 {K}(x,y){u}(0,y)dy =
  (q(x)-p(x))x \www b(x) \\
+ & \int^x_0 K(x,y)(q(x)-p(y))y\www b(y) dy
+ \int^x_0 a(y){K}(x,y)\ppp_y^2(y \www b(y)) dy.
\end{split}
\end{equation}
Consequently, (\ref{saboromond1}) implies
\begin{equation}\label{trosnik1}\begin{split}
&\kappa \vert q(x)-p(x)\vert x\le \vert q(x)-p(x)\vert \vert x \vert 
\www b(x)\vert                                               \\
\le& \int^x_0 \vert K(x,y)(q(x)-p(y))y \www b(y)\vert dy
+ \int^x_0 \vert a(y){K}(x,y)\ppp_y^2(y\www b(y))\vert dy\\
\le & x^2\sup_{y\in [0,x]}\vert K(x,y)\vert x \sup_{y\in [0,x]}\vert q(x)-p(y)
\vert \Vert\www b\Vert_{C[0,\ell]}\\
+ & x\Vert a\Vert_{C[0,\ell]} \sup_{y\in [0,x]}\vert K(x,y)\vert 
\vert \Vert\www b\Vert_{C^2[0,\ell]}.
\end{split}
\end{equation}
The standard a priori estimate for the Goursat problem implies 
\begin{equation}\label{trosnik}
\Vert K\Vert_{C(\overline{\Omega_x})}
\le C\left\Vert \int_0^x (p(y)-q(y))dy \right\Vert_{C[0,x]}
\le C x \Vert p-q\Vert_{C[0,x]}.
\end{equation}
Using (\ref{trosnik}), we can estimate the last terms in (\ref{trosnik1}):
$$
\kappa x\vert p(x)-q(x)\vert \le C x^2 \Vert p-q\Vert_{C[0,x]},
$$
which implies 
$$
\kappa \Vert p-q\Vert_{C[0,x]} \le C x \Vert p-q\Vert_{C[0,x]}.
$$
Hence there exists $\ep>0$ such that 
\begin{equation}\label{saboromond1A} 
p(x)=q(x)\quad \mbox{on}\quad [0,\ep].
\end{equation}
Thus we can always reduces Case A to Case B, so that it suffices to prove 
the theorem only for Case B.
\\
{\bf Proof in Case B.} 
Let $\ep_1\in (0,\ell]$ be the maximal number such that
\begin{equation}\label{xui1} 
p(x)=q(x)\quad \mbox{on}\quad (0,\ep_1).
\end{equation}
Let $\ep_1<\ell$.  Otherwise, the theorem is already proved. 
Let $\delta>0$ be a small parameter. 
We set  $p_\delta(x)=p(x+\ep_1-\delta)$, $q_\delta(x)=q(x+\ep_1-\delta)$,
$a_\delta(x)=a(x+\ep_1-\delta)$, $u_\delta(t,x)=u(t,x+\ep_1-\delta)$
and $\www u_\delta(t,x)=\www u(t,x+\ep_1-\delta)$.

Let a function $\www{v}_\delta $ be given by 
\begin{equation}\label{mkP}
\www{v}_\delta= u_\delta + K_\delta u_\delta.
\end{equation}
Here the operator $K_{\delta}$ is constructed by (\ref{spartak}) 
with the kernel $K$ replaced by the the solution $K_\delta$ 
to the following problem
\begin{equation}\label{(7.1)}
\left\{ \begin{array}{rl}
& a_{\delta}(x)\ppp_x^2{ K}_\delta (x,y) -\ppp_y^2( a_\delta (y)
{ K}_\delta (x,y) )= q_\delta(x){ K}_\delta (x,y) - K_\delta (x,y)
p_\delta (y),\\ 
& \qquad \qquad \quad x\in (0, \ell-\ep_1+\delta),\quad  0 < y < x, \\
& a_\delta(0)\ppp_y {K}_\delta(x,0)+a_\delta'(0){K}_\delta(x,0) = 0, \quad 
      0<x<\ell-\ep_1+\delta, \\
& 2a_\delta (x)\frac{d}{dx}{K}_\delta(x,x)+a_\delta '(x){K}_\delta (x,x) 
= q_\delta(x) - p_\delta (x),\, K_\delta (0,0)=0,\\ 
&\qquad \qquad \quad 0<x<\ell-\ep_1+\delta.
\end{array}\right.   
\end{equation}  
By Lemma 2, we see that $\www{v}_\delta$ satisfies  
\begin{equation}\label{XMmk1}
\begin{split}
&\sigma\ppp_t\www{v}_\delta - a_\delta (x)\ppp_x^2\www{v}_\delta 
+ q_\delta(x)\www{v}_\delta 
= -a_\delta(0) K_\delta(x,0)\partial_x u_\delta(t,0)\quad\mbox{ in}
\quad (0,T)\times (0,\ell-\ep_1+\delta),\\
\quad &\partial_x\www{v}_\delta(t,0)=\partial_x{u}_\delta(t,0)\quad \mbox{and}
\quad  \www{v}_\delta(t,0)=u_\delta(t,0) \quad \mbox{for $0<t<T$}.
\end{split}
\end{equation}

By (\ref{xui1}) and (\ref{(3.1)}), it follows that 
$$
K_\delta(x,x)=0 \quad \mbox{on}\,\, [0,\delta].
$$
Hence this equality and the uniqueness in $\Omega_{\delta}$ of 
solution to (\ref{(3.1)}) imply
\begin{equation}\label{(P3.1)}
\left\{ \begin{array}{rl}
& a_\delta(x)\ppp_x^2{ K}_\delta (x,y) -\ppp_y^2( a_\delta(y){ K}_\delta(x,y))
= q_\delta(x){ K}_\delta(x,y)- { K}_\delta(x,y)p_\delta(y),\\
& \qquad \qquad \quad (x,y) \in \Omega_{\delta}:= \{(x,y);\, 0<y<x<\delta\}, \\
& a_\delta(0)\ppp_y {K}_\delta(x,0)+a_\delta'(0){K_\delta}(x,0) = 0, 
\quad 0<x<\delta, \\
& K_\delta(x,x)=0, \quad 0<x<\delta.
\end{array}\right.
\end{equation}
Therefore, 
\begin{equation}\label{samba}
K_\delta(x,y)=0\quad \mbox{in $\Omega_{\delta}$}.
\end{equation}
Consider the characteristics for the hyperbolic equation in (\ref{(P3.1)}) 
which passes the point $(\delta,\delta)$:
$$
\frac{d\tilde y_\delta}{dx}=-\,\,\root\of{\frac{a(\tilde y_\delta)}{a(x)}},
\quad x\ge \delta,\quad  \tilde y_\delta(\delta)=\delta.
$$
Here we note that we can choose a characteristic curve on which  
$y$ decreases as $x$ increases in a neighborhood of $(\delta, \delta)$.
Then, since $a$ is positive valued, this characteristic curve intersects 
the $x$-axis. By $\widehat\ep(\delta)$, we denote the 
$x$-coordinate of the intersection point.  Then we see that 
$\widehat{\ep}(\delta) > \delta$ and $\lim_{\delta \downarrow 0}
\widehat \ep(\delta) = 0$.
Using the above equality,
from (\ref{(P3.1)}) and (\ref{samba}) we derive
\begin{equation}\label{(KKP3.1)}
\left\{ \begin{array}{rl}
& a_\delta(x)\ppp_x^2{ K}_\delta(x,y) -\ppp_y^2( a_\delta(y){ K}_\delta(x,y))
= q_\delta(x){ K}_\delta(x,y)- { K}_\delta(x,y)p_\delta(y), \\ 
& \qquad \qquad \quad (x,y) \in \{(x,y);\, 0<y<\www y_\delta(x), \,
\delta < x < \widehat\ep(\delta) \}, \\
& a_\delta(0)\ppp_y {K}_\delta(x,0)+a_\delta'(0){K}_\delta(x,0) = 0, 
                  \quad \delta<x<\ep_2(\delta), \\
& K_\delta(\delta,y)=\partial_x K_\delta(\delta,y)=0, \quad 0<y<\delta,
\end{array}\right.
                                        \end{equation}
where we set $\ep_2(\delta) = \mbox{min}\,\{\widehat\ep(\delta),\, \ell\}$.
By the uniqueness of solution for the problem whose initial data are
given on $x=\delta$ in the $xy$-plane, we obtain
$$
K_\delta\equiv 0\quad \mbox{in}\quad \{(x,y);\, 
\delta <y<\widetilde y_\delta(x), \, \delta<x<\ep_2(\delta)\}.
$$
Therefore,
\begin{equation} \label{pharaon} 
K_\delta(x,0)=0\quad \mbox{on}\quad [0,\ep_2(\delta)].
\end{equation} 
In terms of (\ref{pharaon}),  we rewrite the system (\ref{XMmk1}) as 
\begin{equation}\label{PXMmk1}
\begin{split}
&\sigma\ppp_t\www{v}_\delta - a(x)\ppp_x^2\www{v}_\delta 
+ q_\delta(x)\www{v}_\delta =0 \quad\mbox{ in}\quad  
(0,T)\times (0,\ep_2(\delta)),\\
& \partial_x\www{v}_\delta(t,0)=\partial_x{u}_\delta(t,0)\quad 
\mbox{and}\quad \www{v}_\delta(t,0)={u}_\delta(t,0)\quad\mbox{in} 
\,\, (0,T).
\end{split}
\end{equation}  
We apply the unique continuation result which is provided as Lemma \ref{KSOKOL2} in 
Section \ref{SK2} to (\ref{PXMmk1}), so that we can conclude 
$$
\widetilde v_\delta= \widetilde u_\delta \quad \mbox{in}\quad 
(0,T)\times (0,\ep_2(\delta)).
$$
This implies 
$$
K_\delta u_\delta(0,x)=0 \quad \mbox{in}\quad(0,\ep_2(\delta)).
$$

Next we repeat the arguments (\ref{saboromond1})-(\ref{saboromond1A}) 
to obtain
$$
p_\delta=q_\delta\quad \mbox{in}\quad (0,\ep_2(\delta)).
$$
Therefore
$$
p(x)=q(x)\quad\mbox{for $0 < x < \ep_1+\ep_2(\delta)-\delta$}.
$$
Since $\ep_2(\delta)>\delta$, by means of the maximality of $\ep_1$,
we reach a contradiction. 
Thus the proof of the theorem is complete. 
$\blacksquare$

\section{Unique continuation and the proof.}\label{SK2}

In this section, we will prove the following lemma which we used
for the proof of Theorem \ref{Moon} in Section \ref{SNOM}.
We set $Q:= (0,T) \times (0,\ell)$ and 
$H^{1,2}(Q) := \{z;\, z,\partial_t z,\partial_x z,\partial^2_x z\in L^2(Q)\}$.

\begin{lemma}\label{KSOKOL2} 
{\it 
Let $p, q \in C^1[0,\ell], a\in C^2[0,\ell], a(x)>0$ on $[0,\ell]$.
Let $z\in H^{1,2}(Q)$ satisfy 
\begin{equation}\label{kisa}
\left\{ \begin{array}{rl}
& \sigma \ppp_tz - a(x)\ppp_x^2z + p(x)z =0 
\quad \mbox{in}\quad Q,\\
& z(t,0)=\partial_x z(t,0)=0, \quad 0<t<T.
\end{array}\right.
\end{equation}
Then $z\equiv 0$ in $Q$.
}
\end{lemma}

This property called a unique continuation and
can be proved by an $L^2$-weighted estimate called a Carleman estimate.  
We can refer for example to Isakov \cite{Is2} as for general treatments.
However, for our equation \eqref{kisa}, more direct proof is available, 
which is described in this section.

We introduce a weight function $\alpha(t,x)$ of the Carleman estimate:
\begin{equation}\label{DD1.2} 
\alpha(t,x) = e^{\lambda\psi(t,x)}, \quad \psi \in C^2(\ooo Q), 
\quad \vert \partial_x\psi(t,x)\vert\ge C>0\quad \mbox{on}\,\, [0,\ell],
\end{equation}
where $\lambda>0$ is a large parameter.

We define a partial differential operator by
\begin{equation}\label{LD1.4}
P (t,x,D)z := \sigma\ppp_tz -a(x)\ppp_x^2z + p(x)z \quad \mbox{in}\quad Q.
\end{equation}
We introduce the scalar product 
$$
\left\langle u,v\right\rangle_{L^2(Q)} = \int_Q u \ooo{v}dy.
$$
Henceforth, $\ooo{\eta}$ denotes the complex conjugate of $\eta\in \C$.
We are ready to state the main Carleman estimate.
\begin{lemma}\label{SOKOL2} 
{\it Let $p, q \in C^1[0,\ell], a\in C^2[0,\ell], a(x)>0$ on $[0,\ell].$
Let the function $ \alpha$ be defined in \eqref{DD1.2}. 
Then there exists a constant $\widehat\lambda > 0$ 
such that for arbitrary $\lambda \ge \widehat{\lambda}$,
there exists a constant ${\tau}_0(\lambda)>0$ such that 
\begin{equation}\label{DD1.6}
\begin{split}
  \int_Q ( {\tau}|\nabla z|^2 + {\tau}^3\vert z\vert^2)
e^{2{\tau}\alpha}   dx\,dt
\le C_1\int_Q |P(t,x,D)z|^2 e^{2{\tau}\alpha}
dxdt
\end{split}
\end{equation}
for all ${\tau} \ge \tau_0(\lambda)$ and each $z\in H^{1,2}(Q)$ 
satisfying 
\begin{equation}\label{LD1.5}
 z=0\quad \mbox{on}\,\,Q\quad \mbox{and}\quad \partial_x z(t,0)=\partial_x z(t,0)=0\quad \mbox{on}\,\, [0,T].
\end{equation}
}
\end{lemma}

{\bf Proof of Lemma \ref{SOKOL2}.}
We set $g:= P(t,x,D)z$ in $Q$. 
For the proof of the estimate \eqref{DD1.6},
the zeroth order term $p(x)z(t,x)$ of the operator $P$ does not play 
any essential role. 
This term  can be incorporated into the right-hand side, 
and later be absorbed  by the left-hand side of the inequality \eqref{DD1.6}. 
Therefore, it suffices to set
\begin{equation}\label{D1.40}
P(t,x,D)z := \sigma\partial_t z - a(x)\partial^2_x z
\quad \mbox{in $Q$}
\end{equation} 
and consider 
\begin{equation}\label{D1.41}
P(t,x,D)z(t,x) = g(t,x) - p(x)z(t,x) =: \www{g}(t, x) \quad 
\mbox{in $Q$}.
\end{equation}

We denote $w(t,x)=e^{{\tau}\alpha}z(t,x)$.
It follows from \eqref{LD1.4}, 
and \eqref{D1.40}, \eqref{D1.41} that
\begin{equation}\label{D1.44}
P(t,x,D,\tau)w := e^{{\tau}\alpha} P(t,x,D) e^{-{\tau}\alpha} w
= e^{{\tau}\alpha} \www{g} \quad \mbox{in}\quad  Q.
\end{equation}

The operator $P(t,x,D,\tau)$ can be written as follows:
\begin{equation}\label{D1.46}
\begin{split}
&P(t,x,D,\tau)w = \sigma\partial_t w -  a
\partial^2_{x} w 
+ 2 {\tau} a\partial_x\alpha\partial_{x} w                 \\
- & {\tau}^2  a(\partial_x\alpha)^2 w
 + {\tau}  a \partial^2_{x}\alpha w
- {\sigma\tau}(\partial_t\alpha) w.
\end{split}
\end{equation}
We further define two operators by
\begin{equation}\label{sosna}
\left\{ \begin{array}{rl}
& P_1(t,x,D,\tau)u: = \sigma_1\partial_t u+
2\tau a\partial_x\alpha \partial_xu+ {\tau} \partial_x( a \partial_{x}\alpha) w
\cr \\
& P_2(t,x,D,\tau)u=i\sigma_2\partial_t u- \partial_x(a \partial_{x}u )
-  {\tau}^2 a(\partial_x\alpha)^2 u,
\end{array}\right.
\end{equation} 
where $\sigma=\sigma_1+i\sigma_2$ and $\sigma_1,\sigma_2 \in \R$.
From \eqref{D1.46}, we have
\begin{equation}\label{voin}
P_1(t,x,D,\tau) w+P_2(t,x,D,\tau)w=f\quad\mbox{in}\,Q,\quad 
\partial_x w\vert_{\ppp Q}= w\vert_{\partial Q}=0,
\end{equation}
where the function $f$ satisfies
\begin{equation}\label{voinK}
\Vert f\Vert_{L^2(Q)}\le C(\Vert ge^{\tau\alpha}\Vert_{L^2(Q)}+\Vert w\Vert
_{L^2(0,T;H^{1,\tau}(\Omega))}).
\end{equation}
Here we set $\Vert v\Vert_{H^{1,\tau}(\OOO)}:=\root\of{\Vert v\Vert^2_{H^1(\Omega)}+\tau^2\Vert v\Vert^2_{L^2(\Omega)}}$.

Taking the $L^2$-norms of both sides of equation \eqref{voin}, 
we obtain
\begin{equation}\label{suka22}
\Vert P_1w\Vert^2_{L^2(Q)}+\Vert P_2w\Vert^2_{L^2(Q)}+2\mbox{Re}\, 
\left\langle P_1w,P_2w \right\rangle_{L^2(Q)}=\Vert f\Vert^2_{L^2(Q)}.
\end{equation}
Next we compute the third term on the right-hand side of equality 
(\ref{suka22}).  The definition (\ref{sosna}) of the operators  $P_j$ implies
\begin{equation}\label{Ssuka1}
\begin{split}
&\mbox{Re} \left\langle P_1w,P_2w\right\rangle_{L^2(Q)}
= \mbox{Re} \left\langle P_1w,i\sigma_2 \partial_t w\right\rangle_{L^2(Q)}\\
+& \mbox{Re} \left\langle P_1w,-  {\tau}^2 a (\partial_x\alpha)^2 w
\right\rangle_{L^2(Q)}
+\mbox{Re} \left\langle P_1w,- \partial_x(a \partial_{x}w )
\right\rangle_{L^2(Q)}.
\end{split}
\end{equation}
Observe that $P_1^*=-P_1$ and the operator $i\sigma_t\partial_t$ is 
self-adjoint. 
Therefore, using (\ref{LD1.5}), we have
\begin{equation}\label{S1}\begin{split}
&\mbox{Re} \left\langle P_1w,i\sigma_2 \partial_t w\right\rangle_{L^2(Q)}
= \mbox{Re}\left\langle [i\sigma_2\partial_t,P_1]w, w\right\rangle_{L^2(Q)}\\
= &\mbox{Re}\left\langle i\sigma_2 (
2\tau a\partial^2_{tx}\alpha \partial_{x}w+\tau\partial^2_{tx}
(a\partial_x\alpha) w),  w\right\rangle_{L^2(Q)}.
\end{split}
\end{equation}
and 
\begin{equation}\label{S2}\begin{split}
&\mbox{Re} \left\langle P_1w,-  {\tau}^2  a(\partial_x\alpha)^2 w
\right\rangle_{L^2(Q)}                        \\ 
= &- \mbox{Re}\left\langle [{\tau}^2 a(\partial_x\alpha)^2 ,P_1]w, w
\right\rangle_{L^2(Q)}=\left\langle 2\tau^3a\partial_x\alpha\partial_x
(a(\partial_x\alpha)^2) w,w \right\rangle_{L^2(Q)}.
\end{split}
\end{equation}
Here $[P,Q]$ denotes the commutator of two operators: 
$[P,Q] := PQ - QP$.
Short computations imply
\begin{equation}\label{S3}
\begin{split}
&\mbox{Re} \left\langle P_1w,- \partial_x(a \partial_{x}w )
\right\rangle_{L^2(Q)}
= \int_Q \biggl(\frac 12\sigma_1a\partial_t\vert\partial_x w\vert^2
- \tau\partial_x(a^2\partial_x\alpha) \vert \partial_x w\vert^2       \\
+ & 2 {\tau}a\partial_x (a\partial_x\alpha)\vert\partial_{x} w  
\vert^2-{\tau}a\partial_x (a\partial_x\alpha)\vert\partial_{x} w  \vert^2 
+ a\partial_x^2(a\partial_x\alpha) \mbox{Re} \{w\partial_x \bar w\}\biggr)dxdt.
\end{split}
\end{equation}
By (\ref{Ssuka1})-(\ref{S3}), we have 
\begin{equation}\label{KSsuka1}
\begin{split}
&\mbox{Re} \left\langle P_1w,P_2w\right\rangle_{L^2(Q)}\\
= &\mbox{Re}\left\langle i\sigma_2 (
2\tau a\partial^2_{tx}\alpha \partial_{x}w
+\tau\partial^2_{tx}(a\partial_x\alpha) w),  w\right\rangle_{L^2(Q)}
+ \left\langle 2\tau^3a\partial_x\alpha\partial_x(a(\partial_x\alpha)^2) w,w 
\right\rangle_{L^2(Q)}\\
+& \int_Q (a\tau\partial^2_x\alpha \vert \partial_x w\vert^2
+ 2 {\tau}a\partial_x (a\partial_x\alpha)\vert\partial_{x} w  
\vert^2-{\tau}\partial_x (a^2\partial_x\alpha)\vert\partial_{x} w  \vert^2
+ a\partial_x^2(a\partial_x\alpha) \mbox{Re} \{w\partial_x \bar w\}) dxdt.
\end{split}
\end{equation}
Observe that
$$
a\partial_x\alpha\partial_x(a(\partial_x\alpha)^2) \ge C_1(\lambda)>0\quad 
\mbox{on}\quad \ooo{Q}
$$
and
$$
2a\tau\partial^2_x\alpha+a\partial_x (a\partial_x\alpha)
-\partial_x (a^2\partial_x\alpha) \ge C_1(\lambda)>0\quad \mbox{on}\quad 
\ooo{Q}
$$
for all sufficiently large  $\lambda$.

These inequalities and (\ref{KSsuka1}) imply
\begin{equation}\label{KKSsuka1}
\begin{split}
&\mbox{Re} \left\langle P_1w,P_2w\right\rangle_{L^2(Q)}  \\
= &\mbox{Re}\left\langle i\sigma_2 (
2\tau a\partial^2_{tx}\alpha \partial_{x}w+\tau\partial^2_{tx}
(a\partial_x\alpha) w),  w\right\rangle_{L^2(Q)} 
+ C_1\int_Q(\tau \vert \partial_x w\vert^2+\tau^3\vert w\vert^2)dxdt\\
+& \int_Q a\partial_x^2(a\partial_x\alpha) \mbox{Re} \{w\partial_x \bar w\})
dxdt.
\end{split}
\end{equation}
On the other hand, 
\begin{equation}\label{km}\begin{split}
& \left\vert \mbox{Re}\left\langle i\sigma_2
(2\tau a\partial^2_{tx}\alpha \partial_{x}w+\tau\partial^2_{tx}
(a\partial_x\alpha) w),  w\right\rangle_{L^2(Q)}\right\vert 
+ \left\vert \int_Q a\partial_x^2(a\partial_x\alpha) \mbox{Re} 
\{w\partial_x \ooo w\} dxdt\right\vert                 \\ 
\le & C_2(\lambda)\int_Q( \vert \partial_x w\vert^2+\tau^2\vert w\vert^2)dxdt.
\end{split}
\end{equation}
Inequalities (\ref{km}) and (\ref{KKSsuka1}) imply that there exists a 
constant $\tau_0(\lambda)>0$ such that  
\begin{equation}\label{zk}
\mbox{Re} \left\langle P_1w,P_2w\right\rangle_{L^2(Q)}
\ge C_3(\lambda)\int_Q(\tau \vert \partial_x w\vert^2+\tau^3\vert w\vert^2)dxdt
\end{equation}
for all $\tau\ge \tau_0$.
Then inequality follows from (\ref{zk}), (\ref{voinK}) and (\ref{suka22}).
Thus the proof of Lemma \ref{SOKOL2} is complete.
$\blacksquare$

Now, in terms of Lemma \ref{SOKOL2}, we can proceed to 
\\
{\bf Proof of Lemma \ref{KSOKOL2}.} 
Let $(t_0,x_0)\in Q$ belong to $\mbox{supp} \, z$.
 Consider the function $\psi(t,x)=\ell+2 -x-N(t-t_0)^2$
where  $N=\frac{\ell+3}{\mbox{min}\{(T-t_0)^2, t_0^2\}}$.
Hence, 
\begin{equation}\label{L1}
\psi(T,x)<0\quad \mbox{for $0\le x\le \ell$} \quad \mbox{and}\quad  
\psi(0,x)<0\quad \mbox{for $0\le x\le \ell$}.
\end{equation}

Let $\delta\in \left(0, \frac{\ell-x_0}{2}\right)$.
Next we introduce a function by
\begin{equation}\label{L2}
\mu\in C^\infty(\Bbb R),\quad \mu(s)=1\quad \mbox{ for} \quad 
s>\ell+2-x_0-\delta\quad \mbox{and} \quad \mu(s)=0\quad \mbox{ for}\quad  
s<2+\delta.
\end{equation}
We set  $v(t,x)=z(t,x) \mu(\psi(t,x)).$ By (\ref{L1}) and (\ref{L2}), we obtain
$$ 
v(T,\cdot)=v(0,\cdot)=0.
$$

By (\ref{L2}) we have 
$$
\partial_x v(t,\ell)=0, \quad 0<t<T.
$$
Then (\ref{kisa}) yields
\begin{equation}\label{samolet}
Pv=-[\mu(\psi),P]z\quad\mbox{in $Q$}, \quad v\vert_{\partial Q}=0, 
\quad \partial_x v(t,0)=\partial_x v(t,\ell)=0, \, 0<t<T.
\end{equation}
The operator $[\mu(\psi),P]$ is the  first order operator:
$[\mu(\psi),P]=c_1(t,x)\partial_x+c_0(t,x)$, where  
\begin{equation}\label{kisa1}
\mbox{supp}\, c_j\subset \{(t,x);\, \psi(t,x)<2+\delta \},\quad j\in\{0,1\}.
\end{equation}
Applying Lemma \ref{SOKOL2} to equation (\ref{samolet}), we obtain 
\begin{equation}\label{PXDD1.6}
\begin{split}
  \int_Q ( {\tau}|\nabla v|^2 + {\tau}^3\vert v\vert^2)
e^{2{\tau}\alpha}   dx\,dt
\le C_1\int_Q |[\mu(\psi),P]z|^2 e^{2{\tau}\alpha}
dxdt\quad \mbox{for all $\tau\ge\tau_0$}.
\end{split}
\end{equation}
Using (\ref{kisa1}), we estimate 
\begin{equation}\label{P1}
\int_Q |[\mu(\psi),P]z|^2 e^{2{\tau}\alpha}dxdt
\le C_4 e^{2\tau(2+\delta)}\quad \mbox{for all $\tau\ge\tau_0$}.
\end{equation}
We choose a constant $\delta_1>0$ satisfying  
$$
\{ (t,x);\, \vert t-t_0\vert^2 + \vert x-x_0\vert^2 < \delta_1^2\}
\subset  \{(t,x);\, \psi(t,x)> \ell+2-x_0-\delta \}.
$$
Hence, there exists a constant $C_5>0$ such that
\begin{equation}\label{P2}
\int_Q ( {\tau}|\nabla v|^2 + {\tau}^3\vert v\vert^2)e^{2{\tau}\alpha}  dx\,dt
\ge C_5e^{2\tau(\ell+2-x_0-\delta)} \quad \mbox{for all $\tau\ge\tau_0$}.
\end{equation}
By  (\ref{P1}), (\ref{P2}) and (\ref{PXDD1.6}), we have 
$$
C_5e^{2\tau(\ell+2-x_0-\delta)}\le C_4 e^{2\tau(2+\delta)}
\quad \mbox{for all $\tau\ge\tau_0$}.
$$
For large $\tau>0$, we obtain contradiction. Thus the proof of Lemma \ref{KSOKOL2}
is complete. 
$\blacksquare$

{\bf Acknowledgments.}
The work was supported by Grant-in-Aid for Scientific Research (A) 20H00117 
and the Grant-in-Aid for Challenging Research (Pioneering) 21K18142
of the Japan Society for the Promotion of Science.

\end{document}